\documentclass[10pt]{article}
\usepackage{amsfonts}
\usepackage{hyperref}
\usepackage{mathrsfs,enumerate,amssymb,amsmath,color,lineno}
\usepackage{indentfirst}
\usepackage{amsmath}
\usepackage{amssymb}
\usepackage[amsmath,thmmarks]{ntheorem}
\usepackage{graphicx}
\usepackage{tikz}

\newtheorem{definition}{\bf Definition}[section]

\newtheorem{theorem}{\bf Theorem}[section]
\newtheorem{remark}{\bf Remark}[section]
\newtheorem{corollary}{\bf Corollary}[section]
\newtheorem{example}{\bf Example}[section]

\newtheorem{proposition}{\bf Proposition}[section]

\begin{document}
\title{{\textbf{A characterization of a class of border continuous triangular conorms}}}
%\thanks {Supported by the National Natural Science Foundation of China (No.12471440)}}
\author{Zhi-qiang Lai\footnote{\emph{E-mail address}: zhingqlai@foxmail.com}, Xue-ping Wang\footnote{Corresponding author. xpwang1@hotmail.com; fax: +86-28-84761502},\\
\emph{School of Mathematical Sciences, Sichuan Normal University,}\\
\emph{Chengdu 610066, Sichuan, People's Republic of China}}

\newcommand{\pp}[2]{\frac{\partial #1}{\partial #2}}
\date{}
\maketitle
\begin{quote}
{\bf Abstract} Let $T^*:[0,1]^2\rightarrow[0,1]$ be a continuous, non-decreasing and associative function with neutral element, $f: [0,1]\rightarrow [0,1]$ be a strictly monotone function and $f^{(-1)}:[0,1]\rightarrow[0,1]$ be the pseudo-inverse of $f$. This article characterizes the function $T: [0,1]^2 \rightarrow [0,1]$ defined by $T(x,y)=f^{(-1)}(T^*(f(x),f(y)))$ when it is a border continuous triangular conorm.

{\textbf{\emph{Keywords}}:} Triangular conorm; Monotone function; Associative function; Border continuous function \\
\end{quote}

\section{Introduction}
An old problem going back to Abel \cite{NH1826} is to construct two-place real functions from one-place real functions by usual addition operation ``+" (or multiplication ``$\cdot$"), so that it has good algebraic properties, for example, associativity. Abel \cite{NH1826} showed that if $f$ is differentiable and satisfies the system of functional equations
$$f(x,f(y,z))=f(z,f(x,y))=f(y,f(z,x))=f(x,f(z,y))=f(z,f(y,x))=f(y,f(x,z))$$
then there exists a differentiable and invertible function $\psi$ such that
$$\psi(f(x,y))=\psi(x)+\psi(y).$$
The above construction methods of associative functions starting from an one-place function play an indispensable role in theory of solving associative equations \cite{CA2006}. Schweizer and Sklar \cite{BS1961} showed that there exists a continuous, strictly decreasing one-place function $f:(0,1] \to [0,\infty]$ with $f(0^+)=+\infty$ and $f(1)=0$ such that the function $T: [0,1]^2\rightarrow [0,1]$ satisfying $T(x,y)=f^{-1}(f(x)+f(y))$ for all $x,y\in (0,1]$ is a strict triangular norm (t-norm for short)
where $f^{-1}$ is the inverse of $f$ and $f$ is called an additive generator of $T$. Klement, Mesiar and Pap \cite{EP2000} defined an additive generator as a strictly decreasing function $f:[0,1] \to [0,\infty]$ which is right-continuous at $0$ with $f(1)=0$ such that for all $(x,y) \in [0,1]^2$,
$$f(x)+f(y)\in \mbox{Ran}(f)\cup [f(0),\infty]$$
and
\begin{equation}\label{eq1}
T(x,y)=f^{(-1)}(f(x)+f(y))
\end{equation}
where $T$ is t-norm and $f^{(-1)}$ is a pseudo-inverse of $f$, and they discussed the relationship between the continuity of the additive generator and the continuity of corresponding t-norm. Most importantly, they proved that a t-norm is continuous Archimedean if and only if it has a continuous additive generator. The readers can refer to \cite{MC2025,EP1999,AM2004,PV2005,PV2008,PV2010,PV2013,XP2025,YM2024} for further results concerning additive generators of associative functions. For example, Vicen\'{\i}k presented some necessary and sufficient conditions for the generated function given by Eq.\eqref{eq1} being a border continuous t-norm in terms of additive generators \cite{PV2008}, and he further discussed the algebraic and analytical properties of the generated function given by Eq.\eqref{eq1} in terms of additive generators \cite{PV2013}. Following Vicen\'{\i}k's idea, by applying the properties of a strictly monotone function $f:[0,1] \to [0,1]$ Wang and Zhang \cite{XP2025} studied the algebraic and analytical properties of all generated functions $T: [0,1]^2\rightarrow [0,1]$ given by
\begin{equation}\label{eq1.1}
T(x,y)=f^{(-1)}(T^*(f(x),f(y)))
\end{equation}
where $T^*:[0,1]^2\rightarrow [0,1]$ is a continuous, non-decreasing and associative function with neutral element, such as the idempotent elements, the continuity (resp. left-continuity/right-continuity), the associativity and the limit property. Therefore, a natural question arise: can we characterize the function $T$ via Eq.\eqref{eq1.1} when it is a border continuous triangular conorm? This article will answer the question.

The rest of this article is organized as follows. In Section $2$, we present some basic concepts and results. In Sections 3, we give some necessary conditions for the function $T$ given by Eq.\eqref{eq1.1} being a border continuous t-conorm. In Section 4, we characterize the function $T$ defined by Eq.\eqref{eq1.1} when it is a border continuous triangular conorm. A conclusion is drawn in Section 5.

\section{Preliminaries}

In this section, we recall some basic concepts and results.

Let $A$ and $B$ be two sets and define $A\setminus B=\{x\in A| x\notin B\}$.
\begin{definition}[\cite{EP2000}]\label{def2.1}
\emph{A t-norm is a binary operator $T:[0,1]^2\rightarrow [0, 1]$ such that for all $x, y, z\in[0,1]$ the following conditions are satisfied:}

$(T1)$  $T(x,y)=T(y,x)$,

$(T2)$  $T(T(x,y),z)=T(x,T(y,z))$,

\emph{$(T3)$  $T(x,y) \leq T(x,z)$ whenever $y \leq z$,}

$(T4)$  $T(x,1)=x$.

\emph{A binary operator $S:[0, 1]^2 \rightarrow [0, 1]$ satisfying (T1)-(T3) and}

\emph{$(S4)$ $S(x,0)=x$ for all $x\in[0,1]$\\
is said to be a t-conorm. A t-norm $T$ (resp. t-conorm $S$) is said to be border continuous if it is continuous on the boundary of the unit square $[0,1]^2$, i.e., on the set $[0,1]^2\setminus (0,1)^2$. The condition (T4) (resp. (S4)) is called a boundary condition of the t-norm $T$ (resp. t-conorm $S$).}
\end{definition}

\begin{definition}[\cite{EP1999,EP2000}]\label{def2.2}
\emph{Let $a, b, c, d\in [-\infty, \infty]$ with $a<b, c<d$ and $f:[a,b]\rightarrow[c,d]$ be a monotone function. Then the function $f^{(-1)}:[c,d]\rightarrow[a,b]$ defined by
\begin{equation*}
f^{(-1)}(y)=\sup\{x\in [a,b]\mid (f(x)-y)(f(b)-f(a))<0\}\
\end{equation*}
is called a pseudo-inverse of the monotone function $f$.}
\end{definition}

Similarly to the definitions in \cite{EP2000}, we have the following one.
\begin{definition}\label{def2.3}
\emph{Let $S:[0, 1]^2\rightarrow [0, 1]$ be a t-conorm.}

\emph{$(i)$ An element $x\in[0,1]$ is called an idempotent element of $S$ if $S(x,x)=x$. The numbers $0$ and $1$ (which are idempotent elements for each t-conorm $S$) are called trivial idempotent elements of $S$. Each idempotent in $(0,1)$ will be called a non-trivial idempotent element of $S$.}

\emph{$(ii)$ An element $x\in(0,1)$ is called a nilpotent element of $S$ if there exists $n\in N=\{1,2,\cdots,n,\cdots\}$ such that $x_S^{(n)}=1$, where $x_S^{(n)}=S(x_S^{(n-1)},x)$ for $n\geq 2$ and $x_S^{(1)}=x$.}

\emph{$(iii)$ A t-conorm $S$ satisfies the cancellation law if $S(x,y)=S(x,z)$ implies $x=1$ or $y=z$.}

\emph{$(iv)$ A t-conorm $S$ satisfies the conditional cancellation law if $S(x,y)=S(x,z)<1$ implies $y=z$.}

\emph{$(v)$ A t-conorm $S$ is said to be Archimedean if for all $x, y\in(0,1)$ there is an $n\in N$ such that $x_S^{(n)}>y$.}

\emph{$(vi)$ A t-conorm $S$ is called strict if it is continuous and strictly monotone.}

\emph{$(vii)$ A t-conorm $S$ is called nilpotent if it is continuous and each $x\in (0,1)$ is its nilpotent element.}
\end{definition}

\begin{remark}\label{remar2.1}
\emph{It is clear that if a t-conorm $S$ satisfies the conditional cancellation law, then $S(x,y)<S(x,z)$ whenever $y<z$ and $\min\{S(x,y),S(x,z)\}<1$. Moreover, continuous Archimedean t-conorms satisfy the conditional cancellation law.}
\end{remark}

In analogy to Proposition 2.3 in \cite{EP2000}, if $S$ is a continuous t-conorm, then $a \in [0,1]$ is an idempotent element of $S$ if and only if $S(a,x)=\max\{a,x\}$ for all $x \in [0,1]$.

\begin{proposition}[\cite{EP2000}]\label{prop2.1}
  For a function $S:[0,1]^2 \to [0,1]$ the following are equivalent:

  $(i)$ $S$ is a continuous t-conorm;

  $(ii)$ $S$ is uniquely representable as an ordinal sum of continuous Archimedean t-conorms, i.e., there exists a uniquely determined (finite or countably infinite) index set $A$, a family of uniquely determined pairwise disjoint open subintervals $((a_{\alpha},e_{\alpha}))_{\alpha \in A}$ of $[0,1]$ and a family of uniquely determined continuous Archimedean t-conorms $(S_{\alpha})_{\alpha \in A}$ such that $S=(\langle a_{\alpha},e_{\alpha},S_{\alpha}\rangle)_{\alpha \in A}$.
\end{proposition}

A monotone function $f:[0,1] \to [0,1]$ is said to be a generally additive generator of a function $T:[0,1]^2 \to [0,1]$ if Eq.\eqref{eq1.1} holds where $f^{(-1)}$ is the pseudo-inverse of $f$ and $T^*:[0,1]^2 \to [0,1]$ is an associative function with neutral element. Obviously, the function $T$ is always commutative, non-decreasing but need not be associative if $T^*$ is a commutative non-decreasing function.

\subsection{The ranges of generally additive generators}

Let $M\subseteq [0,1]$. Denote the set $[0,1]\setminus  M$ by $M^c$ and write
$$\mathcal{A}=\{M\subseteq[0,1] \mid \mbox{there is a strictly montone function } f:[0,1]\rightarrow[0,1]\mbox{ such that }\mbox{Ran}(f)=M\}.$$

Let $M \in \mathcal{A}$. A pair $(\mathcal{S},C)$ is said to be associated with $M \neq [0,1]$ if $\mathcal{S}=\{[b_k, d_k] \subseteq [0,1]\mid k\in K\}$ is a non-empty countable system of pairwise disjoint intervals of a positive length and $C=\{c_k\in [0,1]\mid k\in K\}$ is a non-empty countable set such that $[b_k,d_k] \cap C=\{c_k\}$ for all $k \in K$ and
\begin{equation*}
M= \{c_k\in [0,1]\mid k\in K\}\cup \left([0,1]\setminus \left(\bigcup_{k\in K}[b_k, d_k] \right)\right).
\end{equation*}
A pair $(\mathcal{S},C)$ is said to be associated with $M=[0,1]$ if $\mathcal{S}=\{[1,1]\}$ and $C=\{1\}$ (see Definition 4 in \cite{PV2005}). We shall briefly write $(\mathcal{S},C)=(\{[b_k, d_k]\mid k\in K\}, \{c_k \mid k\in K\})$ instead of $(\mathcal{S},C)=(\{[b_k, d_k] \subseteq [0, 1]\mid k\in K\}, \{c_k\in [0,1] \mid k\in K\})$. Observe that $(M\setminus C)^c=\cup_{k \in K}[b_k,d_k]$ and $[b_k,d_k] \cap M=[b_k,d_k]\cap C=\{c_k\}$ for all $k\in K$. Further, the set $C$ is always non-empty.

Notice that if $f:[0,1]\rightarrow[0,1]$ is a strictly increasing function with $\mbox{Ran}(f)=M \neq [0,1]$, then the sets $\mathcal{S}=\{[f(x^-),f(x^+)]\mid x\in[0,1], f(x^-)<f(x^+)\}$ and $C=\{f(x)\mid x\in[0,1], f(x^-)<f(x^+)\}$ where $f(x^-)=\lim_{t\to x^-}f(t)$ for all $x\in(0,1]$, $f(0^-)=0$, and $f(x^+)=\lim_{t\to x^+}f(t)$ for all $x\in[0,1)$ and $f(1^+)=1$ have all required properties, i.e., the pair $(\mathcal{S},C)$ is associated with $M \neq [0,1]$.

Let $M\subseteq [0,1]$. A point $x\in [0,1]$ is said to be an accumulation point of $M$ from the left if there exists a strictly increasing sequence $(x_{n})_{n\in N}$ of points $x_{n}\in M$ such that $\lim_{n\rightarrow \infty}x_{n}=x$. A point $x\in [0,1]$ is said to be an accumulation point of $M$ from the right if there exists a strictly decreasing sequence $(x_{n})_{n\in N}$ of points $x_{n}\in M$ such that $\lim_{n\rightarrow \infty}x_{n}=x$. Denote the set of all accumulation points of $M$ from the left (resp. right) by $\mbox{Acc}_{-}(M)$ (resp. $\mbox{Acc}_{+}(M)$). Denote the set of all accumulation points of $M$ from the left and from the right by $\mbox{Acc}_{0}(M)$, i.e., $\mbox{Acc}_{0}(M) = \mbox{Acc}_{-}(M)\cap \mbox{Acc}_{+}(M)$. Denote the set $\mbox{Acc}_{-}(M) \cup \mbox{Acc}_{+}(M)$ by $\mbox{Acc}(M)$.

Let $(\mathcal{S},C)=(\{[b_k, d_k]\mid k\in K\}, \{c_k \mid k\in K\})$ be associated with $M \in \mathcal{A}$. Then the following results are true (see Remark 3 in \cite{PV2005}):

(a) $(M\setminus C) \cap (0,1)= \mbox{Acc}_{0}(M)$,

(b) If $b_k>0$ then $b_k \in \mbox{Acc}_{-}(M)$ and

(c) If $d_k<1$ then $d_k \in \mbox{Acc}_{+}(M)$.

\subsection{An operation $\otimes$ on the range of a generally additive generator}

Let $(\mathcal{S},C)=(\{[b_k, d_k]\mid k\in K\}, \{c_k \mid k\in K\})$ be associated with $M \in \mathcal{A}$. Define a function $F_M:[0,1] \to M$ by
$$\{F_M(x)\}=M \cap [\sup([0,x]\cap M),\inf([x,1]\cap M)] \mbox{ for all } x\in [0,1]$$
(see Definition 5 in \cite{PV2005}). It is a matter of straightforward verification (see Proposition 1 and Lemma 4 in \cite{PV2005}) that
\begin{equation}\label{eq1.2}
F_M(x)=f(f^{(-1)}(x))=\begin{cases}
x & \hbox{if }\ x\in M, \\
c_k &  \hbox{if }\ x\in [b_k,d_k]\setminus \{c_k\} \mbox{ for some } k\in K.
\end{cases}
\end{equation}

\begin{definition}\label{def2.3}
\emph{Let a strictly increasing function $f:[0,1] \to [0,1]$ be a generally additive generator of $T$ and $f^{-1}$ be the inverse of $f$. Define an operation $\otimes: \mbox{Ran}(f)^2 \to \mbox{Ran}(f) $ by
\begin{equation}\label{eq1.3}
x\otimes y=f(T(f^{-1}(x),f^{-1}(y))) \mbox{ for all } x,y \in \mbox{Ran}(f).
\end{equation}}
\end{definition}

Clearly, the operation $\otimes$ is associative if and only if $T$ is associative. Denote $\mbox{Ran}(f)$ by $M$. Substituting Eqs.\eqref{eq1.1} and \eqref{eq1.2} into Eq.\eqref{eq1.3} yields
$$x \otimes y=F_M(T^*(x,y)) \mbox{ for all }x,y \in M$$
where $T^*:[0,1]^2 \to [0,1]$ is an associative function with neutral element. Let $(\mathcal{S},C)=(\{[b_k, d_k]\mid k\in K\}, \{c_k \mid k\in K\})$ be associated with $M \in \mathcal{A}$. Then for all $x,y\in M$ we have that
\begin{equation*}
x \otimes y=\begin{cases}
T^*(x,y) & \hbox{if }\ T^*(x,y)\in M, \\
c_k &  \hbox{if }\ T^*(x,y)\in [b_k,d_k]\setminus \{c_k\} \mbox{ for some } k\in K.
\end{cases}
\end{equation*}

In what follows, the following results are used frequently:

Let $M \in \mathcal{A}$ and $x,y \in M$.

(a) If $T^*(x,y) \leq a$ (resp. $a\leq T^*(x,y))$ and $a\in M$, then $x\otimes y \leq a$ (resp. $a \leq x\otimes y)$.

(b) If $T^*(x,y) < a$ and $a\in \mbox{Acc}_{-}(M)$, then $x\otimes y<a$.

(c) If $a<T^*(x,y)$ and $a\in \mbox{Acc}_{+}(M)$, then $a<x\otimes y$.

\section{The necessary conditions}
This section shows some necessary conditions for the generated function $T$ via Eq.\eqref{eq1.1} being a border continuous t-conorm.

Let $\emptyset \neq A,B\subseteq [0,1]$ and $c\in[0,1]$. Denote $T^*(A,c)=\{T^*(x,c)\mid x\in A\}$, $T^*(c,A)=\{T^*(c,x)\mid x\in A\}$ and $T^*(A,B)=\{T^*(x,y)\mid x\in A, y\in B\}$. In particular, we define $T^*(A,B)=\emptyset$ if $A=\emptyset$ or $B=\emptyset$. By convention, $T^*(A)=A$ for all $A \subseteq [0,1]$.

In what follows, we always suppose that $T^*:[0,1]^2\rightarrow[0,1]$ be a continuous, non-decreasing and associative function with neutral element. Then from Lemma 2.1 of \cite{XP2025} we know that a continuous non-decreasing associative function $T^*:[0,1]^2\rightarrow[0,1]$ with neutral element is either a continuous t-norm or a continuous t-conorm.
\begin{proposition}[\cite{XP2025}]\label{prop3.1}
  Let a function $T$ be given by Eq.\eqref{eq1.1} in which $f:[0,1] \to [0,1]$ is a strictly increasing function and $T^*:[0,1]^2 \to [0,1]$ is a continuous t-conorm. Then the following assertions are equivalent:

  $(i)$ $T$ is continuous at $(0,0)$ and $T(0,0)=0$;

  $(ii)$ $f(0^+)$ is an idempotent element of $T^*$;

  $(iii)$ $T$ is border continuous and satisfies the boundary condition $(S4)$.
\end{proposition}

Note that from Proposition \ref{prop3.1} the main problem for characterizing generally additive generators of border continuous t-conorms becomes to describe the strictly increasing function $f:[0,1] \to [0,1]$ which satisfies that $f(0^+)$ is an idempotent element of $T^*$ and the generated function $T$ via Eq.\eqref{eq1.1} is an associative function. Meanwhile, the associativity of $T$ is equivalent to that of $\otimes$, and the associativity of $\otimes$ depends only on $M$. Therefore, we can reformulate the question mentioned in Section 1 as follows: what is a characterization of $M$ with $M=\mbox{Ran}(f)$, where $f:[0,1] \to [0,1]$ is a strictly increasing function satisfying that $f(0^+)$ is an idempotent element of $T^*$, such that the operation $\otimes$ is associative?

 Let $(\mathcal{S},C)=(\{[b_k, d_k]\mid k\in K\}, \{c_k \mid k\in K\})$ be associated with $M \in \mathcal{A}$. For convenience, if $f(0^+)>0$ then we stipulate $0 \in K, b_0=f(0^-)=0, c_0=f(0)$ and $d_0=f(0^+)$. If $f(0^+)=0$ then we stipulate $0 \notin K$. Write

(D0) $c_k<d_k$ whenever $d_k<1$ with $k\in K \setminus \{0\}$;

(D1) $T^*(M\setminus \{1\},C\setminus \{f(0)\}) \cap (M\setminus C)=\emptyset$;

(D2) $T^*(M\setminus \{1\},\cup_{k \in K\setminus \{0\}} [b_k,c_k]) \cap (M\setminus C)=\emptyset$;

(D3) $T^*(M,M\setminus C)\cap \{d_k|d_k<1,k\in K\setminus \{0\} \}=\emptyset$;

(D4) $(C,\otimes)$ is a semigroup where $\otimes:C \times C \to C$ is the restriction of the operation $\otimes:M \times M \to M$;

(D5) $T^*((M\setminus \{1\}) \cup (\cup_{k \in K\setminus \{0\}} [b_k,c_k]),\cup_{k \in K\setminus \{0\}} [b_k,c_k]) \cap (M\setminus C)=\emptyset$.

\begin{remark}\label{remar3.1}
\emph{It is clear that (D5) implies (D2), and (D2) implies (D1). If $b_k=c_k$ for all $k\in K\setminus \{0\}$ then (D2) coincides with (D1).}
\end{remark}

\begin{proposition}\label{prop3.2}
Let $(\mathcal{S},C)=(\{[b_k, d_k]\mid k\in K\}, \{c_k \mid k\in K\})$ be associated with $M \in \mathcal{A}$, $T^*:[0,1]^2 \to [0,1]$ be a continuous t-conorm and $f:[0,1] \to [0,1]$ be a strictly increasing function with $f(0^+)$ an idempotent element of $T^*$ satisfying that there exist two elements $m_1,n_1 \in [0,1]$ with $m_1<n_1$ such that $f(x) \in (m_1,n_1)$ for each $x\in (0,1)$ and $\langle m_1,n_1,T^*_1 \rangle$ is the summand of the continuous t-conorm $T^*$. If $(M, \otimes)$ is a semigroup, then (D0)-(D5) hold.
\end{proposition}
\begin{proof}
We prove that (D0) is true. Indeed, if there exists an element $k\in K \setminus \{0\}$ such that $c_k=d_k<1$, then $b_k>0$ and $b_k \in \mbox{Acc}_{-}(M)$. From the continuity of $T^*$, there exists an $\varepsilon >f(0^+)$ such that $d_k=T^*(\varepsilon,b_k)$. Because $T^*_1$ is a continuous Archimedean t-conorm, it satisfies the conditionally cancellation law. Thus $b_k=T^*(f(0^+),b_k)<T^*(\varepsilon,b_k)=d_k$ since $f(0^+)$ is an idempotent element of $T^*$. From $f(0^+) \in \mbox{Acc}_{+}(M)$, there exists an $x\in M$ such that $f(0^+)<x \leq \varepsilon$. Hence,  $b_k=T^*(f(0^+),b_k)<T^*(x,b_k) \leq T^*(\varepsilon,b_k)=d_k$, i.e., $b_k<T^*(x,b_k) \leq d_k$. This together with $b_k \in \mbox{Acc}_{-}(M)$ implies that there is a $y \in M$, which is near enough to $b_k$ and satisfies $y<b_k$, such that $b_k \leq T^*(x,y)<T^*(x,b_k) \leq d_k$ since $T^*$ is continuous. Thus $x \otimes y=d_k$ where $\{d_k\}=M \cap [b_k,d_k]$. Again from $f(0^+) \in \mbox{Acc}_{+}(M)$ and the continuity of $T^*$, there exists a $z\in M$, which is near enough to $f(0^+)$ and satisfies $z>f(0^+)$, such that $T^*(y,z)<b_k$ since $y<b_k$. This means $y \otimes z<b_k$ since $b_k \in \mbox{Acc}_{-}(M)$. Therefore, $d_k \geq T^*(x,b_k) \geq T^*(x,y\otimes z)$, which follows that $$d_k \geq x \otimes(y \otimes z).$$ On the other hand, from $z>f(0^+)$ and $T^*_1$ satisfies the conditionally cancellation law, we have $T^*(x \otimes y,z)=T^*(d_k,z)>T^*(d_k,f(0^+))=d_k$ since $f(0^+)$ is an idempotent element of $T^*$. Hence $$(x\otimes y) \otimes z>d_k$$ since $d_k \in \mbox{Acc}_{+}(M)$, a contradiction.

We prove that (D1) is true. If $M=[0,1]$, then $$T^*(M\setminus \{1\},C\setminus \{f(0)\}) \cap (M\setminus C)=T^*([0,1),\{1\}) \cap [0,1)=\{1\} \cap [0,1)=\emptyset.$$ Now, we suppose $M \neq [0,1]$. Assuming there exist $x\in M\setminus \{1\}$ and $c_k \in C\setminus \{f(0)\}$ such that $T^*(x,c_k)\in M\setminus C$. Thus $x\otimes c_k=T^*(x,c_k), k \in K\setminus \{0\}$ and $c_k<1$. Then by (D0), we have $c_k<d_k$ for any $d_k\in [0,1]$. From the continuity of $T^*$, there exists an $\varepsilon >f(0^+)$ such that $d_k=T^*(\varepsilon,c_k)$. Since $f(0^+) \in \mbox{Acc}_{+}(M)$, there is a $z\in M$ such that $f(0^+)<z<\varepsilon$. On the one hand, from the monotonicity of $T^*$, we have $c_k \leq T^*(c_k,z) \leq T^*(c_k,\varepsilon)=d_k$. Thus $c_k \otimes z=c_k$ where $\{c_k\}=M\cap [b_k,d_k]$. Hence $$x \otimes (c_k \otimes z)=x \otimes c_k=T^*(x,c_k).$$ On the other hand, since $T^*_1$ satisfies the conditionally cancellation law and $f(0^+)$ is an idempotent element of $T^*$, we have $T^*(x,c_k)=T^*(T^*(x,c_k),f(0^+))<T^*(T^*(x,c_k),z)=T^*(x \otimes c_k,z)$, i.e., $T^*(x,c_k)<T^*(x \otimes c_k,z)$, which follows that $$T^*(x,c_k)<(x\otimes c_k) \otimes z$$ since $T^*(x,c_k) \in \mbox{Acc}_{+}(M)$, a contradiction.

We prove that (D2) is true. Since (D1) is satisfied, it is enough to prove that $$T^*(M\setminus \{1\},\cup_{k \in K_0} [b_k,c_k)) \cap (M\setminus C)=\emptyset$$ where $K_0=\{k\in K\setminus \{0\}|b_k<c_k\}$. In fact, if there exist $x\in M\setminus \{1\}$ and $k\in K_0$ such that $T^*(x,[b_k,c_k)) \cap (M\setminus C) \neq \emptyset$. From the continuity of $T^*$, we have $[T^*(x,b_k),T^*(x,c_k))\cap (M\setminus C) \neq \emptyset$. Then there is a $p\in [T^*(x,b_k),T^*(x,c_k))\cap (M\setminus C)$, and so $p \in \mbox{Acc}_{+}(M)$. Thus there exist two elements $p_1,p_2\in (p,T^*(x,c_k)) \cap M$ such that $p_1<p_2$. Since $T^*$ is continuous, there is an $\varepsilon>f(0^+)$ such that $c_k=T^*(\varepsilon,b_k)$. From $f(0^+) \in \mbox{Acc}_{+}(M)$ and $p<p_1$, there exists a $z\in M$, which is near enough to $f(0^+)$ and satisfies $f(0^+)<z<\varepsilon$, such that $T^*(p,z)<p_1$. Because $f(0^+)$ is an idempotent element of $T^*$ and $T^*_1$ satisfies the conditionally cancellation law, we have $b_k=T^*(b_k,f(0^+))<T^*(b_k,z)$, i.e., $b_k<T^*(b_k,z)$. From $b_k \in \mbox{Acc}_{-}(M)$, there exists an element $y\in M$, which is near enough to $b_k$ and satisfies $y<b_k$, such that $b_k \leq T^*(y,z)$. On the one hand, $T^*(x,y)<T^*(x,b_k) \leq p$ since $T^*_1$ satisfies the conditionally cancellation law. Hence $x\otimes y \leq p$, which follows that $T^*(x\otimes y,z)\leq T^*(p,z)<p_1$. Therefore, $$(x \otimes y)\otimes z \leq p_1.$$ On the other hand, $b_k \leq T^*(y,z)<T^*(b_k,\varepsilon)=c_k \leq d_k$ since $T^*_1$ satisfies the conditionally cancellation law. Thus $y\otimes z=c_k$ where $\{c_k\}=M \cap [b_k,d_k]$, and so $p_2<T^*(x,c_k)=T^*(x,y\otimes z)$. Hence $$p_2 \leq x\otimes (y\otimes z),$$a contradiction.

We prove that (D3) is true. Suppose that there are $x\in M, y\in M\setminus C$ and $k\in K\setminus \{0\}$ such that $T^*(x,y)=d_k<1$. By (D0), we have $c_k<d_k$.  From the continuity of $T^*$, there exists an element $z>f(0^+)$ such that $d_k=T^*(c_k,z)$. On the one hand, since $f(0^+)$ is an idempotent element of $T^*$ and $T^*_1$ satisfies the conditionally cancellation law, we have $y=T^*(y,f(0^+))<T^*(y,z)$, i.e., $y<T^*(y,z)$. The last inequality means $y<y\otimes z$ since $y \in \mbox{Acc}_{+}(M)$. This implies that $d_k=T^*(x,y)<T^*(x,y \otimes z)$ since $T^*_1$ satisfies the conditionally cancellation law. Hence $$x\otimes(y \otimes z)>d_k$$ since $d_k\in \mbox{Acc}_{+}(M)$. On the other hand, from $T^*(x,y),T^*(c_k,z) \in [b_k,d_k]$, we have $x \otimes y=c_k \otimes z=c_k$ where $\{c_k\}=M\cap [b_k,d_k]$. Thus $(x\otimes y)\otimes z=c_k \otimes z=c_k<d_k$, a contradiction.

We prove that (D4) is true. It suffices to prove that $C\otimes C \subseteq C$. Let $x,y\in C$. Then there are three cases as follows.

(a) If $f(0) \in \{x,y\}$, then there are two subcases as follows.

Subcase 1: If $x=y=f(0)$ then $f(0) \leq T^*(x,y) \leq T^*(f(0^+),f(0^+))=f(0^+)$. Thus $x\otimes y=f(0) \in C$.

Subcase 2: If there is exactly one element of $\{x,y\}$ which is equal to $f(0)$, then $T^*(x,y)=\max\{x,y\} \in C$. Thus $x\otimes y=T^*(x,y) \in C$.

(b) If  $f(0) \notin \{x,y\}$ and $\min\{x,y\}<1$, then $T^*(x,y) \in (M\setminus C)^c=\cup_{k \in K}[b_k,d_k]$ by (D1). Hence there exists an $l\in K$ such that $T^*(x,y) \in [b_l,d_l]$. Thus $x\otimes y=c_l\in C$ where $\{c_l\}=M\cap [b_l,d_l]$.

(c) If  $f(0) \notin \{x,y\}$ and $\min\{x,y\}=1$, then $x\otimes y=1 \in C$.

Now, we prove that (D5) and (D2) are equivalent. From Remark \ref{remar3.1} we just show that (D2) implies (D5). Suppose that (D2) holds. If $M=[0,1]$, then (D5) is true obviously. Let $M\neq [0,1]$. It is enough to prove $[T^*(b_i,b_j),T^*(c_i,c_j)]\cap (M\setminus C)=\emptyset$ for all $i,j \in K\setminus \{0\}$. We assert that \begin{equation*}[T^*(b_i,c_j),T^*(c_i,c_j)] \cap (M\setminus C)=\emptyset\end{equation*} for all $i,j \in K\setminus \{0\}$. Indeed, for any $i,j \in K\setminus \{0\}$, if $c_j=1$ then $[T^*(b_i,c_j),T^*(c_i,c_j)] \cap (M\setminus C)=\emptyset$. If $c_j<1$ then $c_j \in M \setminus \{1\}$, thus from (D2) we have $T^*(c_j,[b_i,c_i]) \cap (M\setminus C)=\emptyset$. This together with the continuity and commutativity of $T^*$ yields that $[T^*(b_i,c_j),T^*(c_i,c_j)]\cap (M\setminus C)=\emptyset$.

It remains to prove $[T^*(b_i,b_j),T^*(b_i,c_j))\cap (M\setminus C)=\emptyset$ for all $i,j \in K\setminus \{0\}$. Obviously, it is true when $b_j=c_j$. Now, suppose that $b_j<c_j$. Since $i \in K\setminus \{0\}$, $b_i \in \mbox{Acc}_{-}(M)$. Choose a strictly increasing sequence $(x_n)_{n \in N}$ in $M\setminus\{1\}$ such that $\lim_{n\to \infty}x_n=b_i$. Thus from (D2), we have $T^*(x_n,[b_j,c_j])\cap (M\setminus C)=\emptyset$ for all $n\in N$. It follows from the monotonicity and continuity of $T^*$ that $[T^*(x_n,b_j),T^*(x_n,c_j)]\cap (M\setminus C)=\emptyset$ for all $n\in N$. Thus $\cup_{n\in N} [T^*(x_n,b_j),T^*(x_n,c_j)]\cap (M\setminus C)=\emptyset$. Since $[T^*(b_i,b_j),T^*(b_i,c_j)) \subseteq \cup_{n\in N} [T^*(x_n,b_j),T^*(x_n,c_j)]$, we have $[T^*(b_i,b_j),T^*(b_i,c_j))\cap (M\setminus C)=\emptyset$. This completes the proof.
\end{proof}

\begin{remark}\label{remar3.2}
\emph{From the proof of Proposition \ref{prop3.2}, the equivalence between (D2) and (D5) does not depend on the associativity of $\otimes$.}
\end{remark}

 The following corollary is a straightforward consequence of Propositions \ref{prop3.1} and \ref{prop3.2}

\begin{corollary}\label{coro3.1}
Let $T^*:[0,1]^2 \to [0,1]$ be a continuous t-conorm and $f:[0,1] \to [0,1]$ be a strictly increasing function satisfying that there exist two elements $m_1,n_1 \in [0,1]$ with $m_1<n_1$ such that $f(x) \in (m_1,n_1)$ for each $x\in (0,1)$ and $\langle m_1,n_1,T^*_1 \rangle$ is the summand of the continuous t-conorm $T^*$. If the generated function $T$ via Eq.\eqref{eq1.1} is a border continuous t-conorm, then $f(0^+)$ is an idempotent element of $T^*$ and (D0)-(D5) hold.
\end{corollary}

In Proposition \ref{prop3.2} and Corollary \ref{coro3.1}, the condition that there exist two elements $m_1,n_1 \in [0,1]$ with $m_1<n_1$ such that $f(x) \in (m_1,n_1)$ for each $x\in (0,1)$ and $\langle m_1,n_1,T^*_1 \rangle$ is the summand of the continuous t-conorm $T^*$ can't be deleted generally.

\begin{example}\label{exp3.1}
\emph{Let the function $S_p:[0,1]^2\rightarrow [0,1]$ be defined by
\begin{equation*}
S_p(x,y)=x+y-xy,
\end{equation*}}
\emph{the function $T^*:[0,1]^2\rightarrow [0,1]$ be defined by
\begin{equation*}
T^*(x,y)=\begin{cases}
\frac{1}{2}S_p(\frac{x-0}{\frac{1}{2}-0},\frac{y-0}{\frac{1}{2}-0}) & \text{if } (x,y)\in [0,\frac{1}{2}]^2,\\
(1-\frac{1}{2})S_p(\frac{x-\frac{1}{2}}{1-\frac{1}{2}},\frac{y-\frac{1}{2}}{1-\frac{1}{2}})+\frac{1}{2} & \text{if } (x,y)\in [\frac{1}{2},1]^2,\\
\max\{x,y\} & \text{otherwise}
\end{cases}
\end{equation*}}
\emph{and the function $f:[0,1] \rightarrow [0,1]$ be defined by}
\emph{\begin{equation*}
f(x)=\begin{cases}
\frac{4}{5}x & \text{if } x\in [0,\frac{1}{2}],\\
\frac{1}{3}(2x+1) & \text{if } x\in (\frac{1}{2},1],
\end{cases}
\end{equation*}
respectively.}

\emph{It is easy to check that $M=[0,\frac{2}{5}] \cup (\frac{2}{3},1]$, $b_k=c_k=\frac{2}{5},d_k=\frac{2}{3}$ and $f(0^+)=0$ is an idempotent element of $T^*$. Meanwhile, $T^*=(\langle0,\frac{1}{2},S_p\rangle,\langle \frac{1}{2},1,S_p\rangle)$ is a continuous t-conorm, $f(x) \in [0,\frac{1}{2}]$ for each $x \in [0,\frac{1}{2}]$, $f(x) \in [\frac{1}{2},1]$ for each $x\in (\frac{1}{2},1]$ and
\begin{equation*}
T(x,y)=f^{(-1)}(T^*(f(x),f(y)))=\begin{cases}
x+y-\frac{8}{5}xy & \text{if } (x,y)\in [0,\frac{1}{2}]^2 \mbox{ and } 0 \leq x+y-\frac{8}{5}xy \leq \frac{1}{2},\\
\frac{1}{2} & \text{if } (x,y)\in [0,\frac{1}{2}]^2 \mbox{ and } \frac{1}{2}<x+y-\frac{8}{5}xy \leq \frac{3}{5},\\
\frac{1}{3}(4x+4y-4xy-1) & \text{if } (x,y)\in (\frac{1}{2},1]^2,\\
\max\{x,y\} & \text{otherwise}
\end{cases}
\end{equation*}
is a border continuous t-conorm, i.e., $(M,\otimes)$ is a semigroup. However, by a simple calculation, we have:}

\emph{$(i)$ $T^*(M\setminus \{1\},C\setminus \{f(0)\})=[\frac{2}{5},\frac{12}{25}]\cup (\frac{2}{3},1)$, $M\setminus C=[0,\frac{2}{5}) \cup (\frac{2}{3},1]$. Hence $T^*(M\setminus \{1\},C\setminus \{f(0)\}) \cap (M\setminus C)=(\frac{2}{3},1) \neq \emptyset$, contrary to (D1).}

\emph{$(ii)$ $T^*(M\setminus \{1\},\cup_{k \in K\setminus \{0\}} [b_k,c_k])=[\frac{2}{5},\frac{12}{25}]\cup (\frac{2}{3},1)$. Hence $T^*(M\setminus \{1\},\cup_{k \in K\setminus \{0\}} [b_k,c_k]) \cap (M\setminus C)=(\frac{2}{3},1) \neq \emptyset$, contrary to (D2) and (D5).}
\end{example}

\section{The necessary and sufficient conditions}
 This section shows a necessary and sufficient condition for the generated function $T$ via Eq.\eqref{eq1.1} being a border continuous t-conorm.

Let $ \mathcal{A}_0=\{M \in \mathcal{A} \mid f(0^+) \notin \mbox{Acc}_{+}(M^c) \}$.
It is clear that $M \in \mathcal{A}_0$ if and only if there exists an $\varepsilon \in (f(0^+),1)$ such that $(f(0^+),\varepsilon] \subseteq M$. We first give the following definition.

 \begin{definition}\label{def4.1}
 \emph{Let $M \in \mathcal{A}_0$. The set $M$ is said to be determined by a sequence of triples $\{(a_i,b_i,c_i)\}_{i\in I}$ if
 $$ M=\{f(0)\} \cup \left(\bigcup_{i\in I}((a_i,b_i)\cup \{c_i\})\right), I=\{1,2,\dotsi,n\}\mbox{ with }n\in N, \mbox{ and } b_n<1 \mbox{ whenever } n \geq 2$$
 or
 $$ M=\{f(0)\} \cup \left(\bigcup_{i\in I}((a_i,b_i)\cup \{c_i\})\right)\cup\{1\},I=N \mbox{ and } \lim_{i\to \infty}c_i=1,$$
 where $(a_i)_{i\in I},(b_i)_{i\in I}$ and $(c_i)_{i\in I}$ are strictly increasing sequences of non-negative real numbers of $[0,1]$ with $a_1=f(0^+), a_i<b_i\leqslant c_i$ for all $i \in I$ and $b_i<a_{i+1}, c_i\leqslant a_{i+1}$ for all $i,i+1 \in I$, respectively,}
 \end{definition}
\begin{example}
\emph{$(1)$ The set $M=\{\frac{1}{100}\} \cup (\frac{1}{4},\frac{1}{2}] \cup (\frac{2}{3},\frac{5}{6}] \cup (\frac{7}{8},\frac{8}{9}]$ is determined by the sequence of triples $\{(\frac{1}{4},\frac{1}{2},\frac{1}{2}),(\frac{2}{3},\frac{5}{6},\frac{5}{6}),(\frac{7}{8},\frac{8}{9},\frac{8}{9})\}$.}

\emph{$(2)$ The set $M=\{0\} \cup \left(\bigcup_{n \in N} (1-\frac{1}{2^n}-\frac{1}{3^n},1-\frac{1}{2^n}]\right) \cup \{1\}$ is determined by the sequence of triples $\{(1-\frac{1}{2^i}-\frac{1}{3^i},1-\frac{1}{2^i},1-\frac{1}{2^i})\}_{i \in N}$.}

\emph{$(3)$ The set $M=[0,\frac{1}{3}]\cup (\frac{1}{2},1] \in \mathcal{A}_0$ but it is not in the form described by Definition \ref{def4.1}.}
\end{example}

 The following are divided into two cases $f(0^+)>0$ and $f(0^+)=0$, respectively.
 \begin{proposition}\label{prop4.1}
 Let $(\mathcal{S}',C')=(\{[b_k', d_k']\mid k\in K\}, \{c_k' \mid k\in K\})$ be associated with $M \in \mathcal{A}_0$, $T^*:[0,1]^2 \to [0,1]$ be a continuous t-conorm, and $f:[0,1] \to [0,1]$ be a strictly increasing function with $f(0^+)>0$ being an idempotent element of $T^*$. If there exist two elements $m_1,n_1 \in [0,1]$ with $m_1<n_1$ such that $f(x) \in (m_1,n_1)$ for each $x\in (0,1)$ and $\langle m_1,n_1,T^*_1 \rangle$ is the summand of the continuous t-conorm $T^*$, then (D2)-(D4) imply the following statements:

 (D6) There exist uniquely determined strictly increasing sequences $(a_i)_{i\in I},(b_i)_{i\in I}$, and $(c_i)_{i\in I}$ of non-negative real numbers of $[0,1]$ such that $M$ is determined by the sequence of triples $\{(a_i,b_i,c_i)\}_{i\in I}$.

 (D7) Either $a_k\leq T^*(a_i,a_j)$ and $b_k\leq \min\{T^*(a_i,b_j),T^*(a_j,b_i)\}$ or $T^*(c_i,c_j)\leq a_k$ for all $i,j,k \in I$.

 (D8) $(\{c_i \mid i\in I\cup \{0\}\}, \otimes)$ is a semigroup where $\otimes: \{c_i \mid i\in I\cup \{0\}\} \times \{c_i \mid i\in I\cup \{0\}\} \to \{c_i \mid i\in I\cup \{0\}\}$ is the restriction of the operation $\otimes: M \times M \to M$.
  \end{proposition}
\begin{proof}
	We first prove that (D6) is true. From Remark \ref{remar3.1}, (D2) implies (D1).
Since $f(0^{+})\notin \mbox{Acc}_{+}(M^{c})$, there exists an $\varepsilon \in (f(0^+),1)$ such that $(f(0^{+}),\varepsilon]\subseteq M$. This together with (D1) implies $T^{*}((f(0^{+}),\varepsilon], c_k')\cap(M \setminus C')=\emptyset$ for all $k\in K \setminus \{0\}$. From the continuity of $T^*$ and $f(0^{+})$ is an idempotent element of $T^{*}$, we have $T^{*}((f(0^{+}),\varepsilon], c_{k}')\cap (M\setminus C')=(T^{*}(f(0^{+}),c_{k}'),T^{*}(\varepsilon,c_{k}')]\cap(M \setminus C')=(c_{k}',T^{*}(\varepsilon,c_{k}')]\cap(M\setminus C')=\emptyset.$ Thus there exists an $l\in K\setminus \{0\}$ such that $(c_{k}',T^{*}(\varepsilon,c_{k}')] \subseteq [b_{l}',d_{l}']$, and consequently, $[c_{k}',T^{*}(\varepsilon,c_{k}')] \subseteq [b_{l}',d_{l}']$. Since $c_k' \in [b_k',d_k']\cap [b_l',d_l']$ and intervals of $\mathcal{S}'$ are pairwise disjoint, we have that $k=l$ and $T^*(\varepsilon,c_k')\le d_k'$. Hence, for all $k_1,k_2\in K \setminus \{0\}$ with $k_1<k_2$, we have $T^*(\varepsilon,c_{k_{1}}')\le d_{k_{1}}' \le c_{k_{2}}'$, i.e., $T^{*}(\varepsilon,c_{k_{1}}')\le c_{k_{2}}'$.

We assert $\mbox{Acc}(C')\subseteq \{1\}$. Indeed, if $a \in \mbox{Acc}(C')$ and $a<1$, then $a=\sup\{c_{k}' \mid k\in K \setminus \{0\}\}$ and $a\in \mbox{Acc}_{-}(C')$. Since $f(0^{+})$ is an idempotent element of $T^{*}$ and $T_{1}^*$ satisfies the conditional cancellation law, we have $a=T^{*}(f(0^+),a)<T^{*}(\varepsilon,a)$. From $a\in \mbox{Acc}_{-}(C')$ and the continuity of $T^*$, there exists a $k_{1}\in K\setminus \{0\}$ such that $c_{k_{1}}'$ is near enough to $a$, $c_{k_{1}}'<a$ and $a<T^*(\varepsilon,c_{k_{1}}')$. Choose a $k_{2}\in K \setminus\{0\}$ with $k_{1}<k_{2}$. Then we have $ c_{k_{1}}'<c_{k_{2}}' \le a$. Thus $c_{k_2}'\leq a< T^*(\varepsilon,c_{k_{1}}')$, a contradiction. Therefore, $\mbox{Acc}(C')\subseteq \{1\}$.

Note that the set $C'$ is always non-empty and countable. Hence, we consider two cases as below.

(a) If $\mbox{Acc}(C')=\emptyset$, then the set $C'$ is finite and can be expressed in the form $C'=\{c_i \mid i\in I \cup \{0\} \}$ where $I=\{1,2,\cdots,n\}$ with $n \in N, c_0=c'_0=f(0)$ and $c_i<c_{i+1}$ for all $i,i+1 \in I \cup \{0\}$.

(b) If $\mbox{Acc}(C')=\{1\}$, then $C'$ is infinite and can be expressed in the form $C'=\{c_i \mid i\in I \cup \{0\} \}$ where $I=N, c_0=c'_0=f(0)$ and $c_i<c_{i+1}$ for all $i,i+1 \in I \cup \{0\}$.

In both cases there exists a uniquely $k\in K$ such that $c_i=c_k' \in[b_{k}',d_{k}']$ for each $i\in I \cup \{0\}$. Put $b_{i}=b_{k}'$ and $d_{i}=d_{k}'$. Then obviously we can get that $C'=\{c_{i} \mid i\in I\cup \{0\}\}$ and $\mathcal{S}'=\{[b_{i},d_{i}] \mid i\in I\cup \{0\}\}$. Finally, put $a_1=f(0^+)$ and $a_{i+1}=d_i$ with $ d_i<1$ for all $i \in I$.

 In case $I=\{1,2,\cdots,n\}$ with $n \in N$ we have $b_n<d_n \leq 1$ whenever $n \geq 2$.

 In case $I=N$, we have proved that $\lim_{i \to \infty}c_i=1$.

Thus, by Definition \ref{def4.1} we know that the set $M$ is determined by a sequence of triples $\{(a_i,b_i,c_i)\}_{i\in I}$.

Now, suppose that the set $M$ is  determined by sequences of triples $\{(a_{i},b_{i},c_{i})\}_{i\in I}$ and $\{(u_{j},v_{j},w_{j})\}_{j\in J}$, respectively. Then clearly, $I=J$ and $a_i=u_j, b_i= v_j$ and $c_i=w_j$ by induction showing the uniqueness of sequences $(a_i)_{i\in I},(b_i)_{i\in I}$ and $(c_i)_{i\in I}$. Therefore, (D6) is true.

Because of $C'=\{c_i \mid i\in I \cup \{0\} \}$, we have that (D4) coincides with (D8). Therefore, (D8) is true.

Next, we prove (D7) is true. Since (D2) is equivalent to (D5) by Remark \ref{remar3.2},  (D5) and (D3) can be equivalently written, respectively, as follows:

\begin{equation}\label{eq4.1}
T^*((M \setminus \{1\}) \cup (\cup_{i\in I}[b_i,c_i]),\cup_{i\in I}[b_i,c_i]) \cap (M\setminus \{c_i \mid i\in I \cup \{0\}\})= \emptyset
\end{equation}
and
\begin{equation}\label{eq4.2}
T^*(M,M\setminus \{c_i \mid i\in I \cup \{0\}\}) \cap \{a_i \mid i \in I\setminus \{1\}\} =\emptyset.
\end{equation}
Substitute $M=\{f(0)\} \cup (\cup_{i\in I}((a_i,b_i)\cup \{c_i\}))$ whenever $I=\{1,2,\dotsi,n\}$ with $n\in N$, and $M=\{f(0)\} \cup (\cup_{i\in I}((a_i,b_i)\cup \{c_i\})) \cup \{1\}$ whenever $I=N$ into Eqs.\eqref{eq4.1} and \eqref{eq4.2}, respectively. Then from the continuous of $T^*$, $T_1^*$ satisfies conditional cancellation law and $f(0^+)$ is an idempotent element of $T^*$ we know that  Eqs.\eqref{eq4.1} and \eqref{eq4.2} imply
\begin{equation}\label{eq4.3}
(\min\{T^*(a_i,b_j),T^*(a_j,b_i)\},T^*(c_i,c_j)) \cap (a_k,b_k)=\emptyset \mbox{ for all } i,j,k\in I
\end{equation}
and
\begin{equation}\label{eq4.4}
a_k \notin (T^*(a_i,a_j),T^*(b_i,b_j)) \cup (T^*(a_i,c_j),T^*(b_i,c_j)) \mbox{ for all } i,j,k\in I,
\end{equation}
respectively.

Below, we show that Eqs.\eqref{eq4.3} and \eqref{eq4.4} are equivalent to
\begin{equation}\label{eq4.5}
(\min\{T^*(a_i,b_j),T^*(a_j,b_i)\},T^*(c_i,c_j)) \cap (a_k,b_k)=\emptyset \mbox{ for all } i,j,k\in I
\end{equation}
and
\begin{equation}\label{eq4.6}
a_k \notin (T^*(a_i,a_j),T^*(b_i,b_j)) \mbox{ for all } i,j,k\in I,
\end{equation}
respectively.
It suffices to prove that Eqs.\eqref{eq4.5} and \eqref{eq4.6} imply that $a_k \notin (T^*(a_i,c_j),T^*(b_i,c_j))$ for all $i,j,k\in I$. In fact, by Eq.\eqref{eq4.5} we know that $a_k \notin (T^*(a_i,b_j),T^*(c_i,c_j))$. This together with $(T^*(a_i,c_j),T^*(b_i,c_j)) \subseteq (T^*(a_i,b_j),T^*(c_i,c_j))$ yields that $a_k \notin (T^*(a_i,c_j),T^*(b_i,c_j))$. Thus Eqs.\eqref{eq4.3} and \eqref{eq4.4} are equivalent to Eqs.\eqref{eq4.5} and \eqref{eq4.6}, respectively.

 Because Eqs.\eqref{eq4.5} and \eqref{eq4.6} are equivalent to (D7), (D7) is true.
\end{proof}

Sequently we come to our main results.
\begin{proposition}\label{propth4.1}
Let $M \in \mathcal{A}_0$, $T^*:[0,1]^2 \to [0,1]$ be a continuous t-conorm, and $f:[0,1] \to [0,1]$ be a strictly increasing function with $f(0^+)>0$ being an idempotent element of $T^*$. If there exist two elements $m_1,n_1 \in [0,1]$ with $m_1<n_1$ such that $f(x) \in (m_1,n_1)$ for each $x\in (0,1)$ and $\langle m_1,n_1,T^*_1 \rangle$ is the summand of the continuous t-conorm $T^*$, then $(M, \otimes)$ is a semigroup if and only if (D6)-(D8) hold.
\end{proposition}
\begin{proof}
$(\Rightarrow)$ If $(M, \otimes)$ is a semigroup then by Propositions \ref{prop3.2} and \ref{prop4.1}, (D6)-(D8) hold.

$(\Leftarrow)$ Supposing that (D6)-(D8) are satisfied, we prove that $(M, \otimes)$ is a semigroup. For all $i\in I$, put $d_i=a_{i+1}$ if $i+1 \in I$ and $d_i=1$ if $i+1 \notin I$. Then  $(\mathcal{S},C)=(\{[b_i, d_i]\mid i\in I\cup \{0\}\}, \{c_i \mid i\in I\cup \{0\}\})$ is associated with $M \in \mathcal{A}_0$ where $b_0=f(0^-)=0,c_0=f(0)$ and $d_0=f(0^+)$. Write

$$m_{i,j}=\min\{T^*(a_i,b_j),T^*(a_j,b_i)\} \mbox{ for all } i,j\in I$$
and define an operation $\oplus:I \times I \to I$ by
\begin{equation}\label{eq4.7}
i \oplus j=\max\{k \in I \mid a_k \leq T^*(a_i,a_j)\}.
\end{equation}
We prove the following three results.

$(i)$ For arbitrary $i,j \in I$,
\begin{equation}\label{eq4.8}
F_{M}(m_{i,j})=c_{i\oplus j}
\end{equation}
and
\begin{equation}\label{eq4.9}
c_{i\oplus j}=c_i\otimes c_j.
\end{equation}

Indeed, from Eq.\eqref{eq4.7} we have $a_{i\oplus j+1}>T^*(a_i,a_j)$ for arbitrary $i,j \in I$. By (D7) we know that $T^{*}(c_{i},c_{j}) \leq a_{i\oplus j+1}=d_{i\oplus j}$, and $b_{i \oplus j} \le m_{i,j} \le T^{*}(c_i,c_j) \le d_{i \oplus j}$ since $T^*$ is monotone. Hence $F_{M}(m_{i,j})=c_{i \oplus j}$ and $c_i \otimes c_j=F_{M}(T^*(c_i,c_j))=c_{i\oplus j}$ where $\{c_{i \oplus j}\}=M \cap [b_{i \oplus j},d_{i \oplus j}]$.

$(ii)$ For arbitrary $x\in (a_i,b_i)\cup \{c_i\},y=c_j$ with $i,j\in I$,
\begin{equation}\label{eq4.10}
x \otimes y=c_i \otimes c_j.
\end{equation}

Indeed, for any $x\in (a_i,b_i)\cup \{c_i\}$ and $y=c_j$ we have $m_{i,j} \leq T^*(a_i,b_j) \leq T^*(x,y) \leq T^*(c_i,c_j)$ since $T^*$ is monotone. Consequently, from Eqs.\eqref{eq4.8} and \eqref{eq4.9} we have $c_{i\oplus j} = F_{M}(m_{i,j}) \leq F_{M}(T^*(x,y)) = x \otimes y \leq F_{M}(T^*(c_i,c_j))=c_i \otimes c_j=c_{i \oplus j}$. Therefore, $x \otimes y=c_{i\oplus j}=c_i \otimes c_j$.

$(iii)$ For arbitrary $x \in (a_i,b_i)\cup \{c_i\}, y\in (a_j,b_j) \cup \{c_j\}$ with $i,j\in I$ satisfying $T^*(x,y) \notin \cup_{l \in I}(a_l,b_l)$,
\begin{equation}\label{eq4.11}
b_{i\oplus j} \leq T^*(x,y)
\end{equation}
and
\begin{equation}\label{eq4.12}
x \otimes y=c_i \otimes c_j.
\end{equation}

Indeed, for arbitrary $x \in (a_i,b_i)\cup \{c_i\}, y\in (a_j,b_j) \cup \{c_j\}$ satisfying $T^*(x,y) \notin \cup_{l \in I}(a_l,b_l)$, we have $a_{i\oplus j} \leq T^*(a_i,a_j)< T^*(x,y) \leq
T^*(c_i,c_j) \leq a_{i\oplus j+1}=d_{i\oplus j}$ since (D7) holds and $T_1^*$ satisfies the conditional cancellation law, which together with $T^*(x,y) \notin \cup_{l \in I}(a_l,b_l)$ implies $T^{*}(x,y) \notin (a_{i\oplus j},b_{i\oplus j})$. Thus $T^*(x,y) \in [b_{i\oplus j},d_{i\oplus j}]$, i.e, $b_{i\oplus j} \leq T^*(x,y)$, and consequently $x \otimes y=c_{i \oplus j}=c_i \otimes c_j$ where $\{c_{i \oplus j}\}=M \cap [b_{i \oplus j},d_{i \oplus j}]$.

Now, we prove that $(M, \otimes)$ is a semigroup. For arbitrary $x,y,z \in M$, there are three cases as follows.

$(a)$ If at least one element of $\{x,y,z\}$ is $f(0)$ or 1, then $(x\otimes y)\otimes z=x \otimes (y\otimes z)$.

$(b)$ If $x,y,z \in M \cap (f(0),1)$ satisfying that $T^*(x,y),T^*(y,z) \in \cup_{l\in I}(a_{l},b_{l})$, then $x\otimes y=T^{*}(x,y)$ and $y\otimes z=T^{*}(y,z)$. On the one hand, $(x \otimes y)\otimes z= F_{M}(T^*(x \otimes y,z)=F_{M}(T^*(T^*(x,y),z))$. On the other hand, $x \otimes (y\otimes z)=F_{M}(T^*(x,y\otimes z))=F_{M}(T^*(x,T^{*}(y,z)))$. Hence $(x\otimes y)\otimes z=x \otimes (y\otimes z)$ since $T^*$ is associative.

$(c)$ If $x,y,z \in M \cap (f(0),1)$ satisfying that $T^*(x,y) \notin \cup_{l\in I}(a_{l},b_{l})$ or $T^*(y,z) \notin \cup_{l\in I}(a_{l},b_{l})$, then there exist three elements $i,j,k \in I$ such that $x \in (a_i,b_i)\cup \{c_i\}, y\in (a_j,b_j)\cup \{c_j\}$ and $z \in(a_k,b_k)\cup \{c_k\}$. Consider three subcases as below,

Subcase 1: If $T^*(x,y),T^*(y,z) \notin \cup_{l\in I}(a_{l},b_{l})$, then from Eq.\eqref{eq4.12}, $x \otimes y=c_i\otimes c_j$. Thus from Eq.\eqref{eq4.9}, we have $c_i\otimes c_j=c_{i\oplus j} \in C$. Hence, by Eq.\eqref{eq4.10} we know that $c_{i\oplus j} \otimes z=c_{i\oplus j}\otimes c_k=(c_i\otimes c_j)\otimes c_k$. This follows that $(x\otimes y)\otimes z=(c_i\otimes c_j) \otimes z=c_{i\oplus j} \otimes z=(c_i\otimes c_j)\otimes c_k$. Similarly, we have $x\otimes (y\otimes z)=c_i\otimes (c_j \otimes c_k)$. Therefore, from (D8) we have that $(x\otimes y)\otimes z=x \otimes (y\otimes z)$.

Subcase 2: If $T^*(x,y) \notin \cup_{l\in I}(a_{l},b_{l}), T^*(y,z) \in \cup_{l\in I}(a_{l},b_{l})$, then similarly to Subcase 1 we have
\begin{equation}\label{eq4.13}
(x\otimes y)\otimes z=(c_i\otimes c_j)\otimes c_k.
\end{equation}
From the inequalities $x \leq c_i$, $y \leq c_j$ and $z\leq c_k$ and the monotonicity of $\otimes$, we have $x\otimes (y\otimes z)\leq c_i \otimes (c_j \otimes c_k)$. This together with (D8) and Eq.\eqref{eq4.13} implies that
$$x\otimes (y\otimes z) \leq c_i \otimes (c_j \otimes c_k)=(c_i \otimes c_j) \otimes c_k=(x\otimes y)\otimes z.$$
Thus it remains to prove $(x\otimes y)\otimes z \leq x\otimes (y\otimes z)$. By $a_k<z$ and \eqref{eq4.11}, we know that $m_{i\oplus j,k} \leq T^*(b_{i \oplus j},a_k) \leq T^*(T^*(x,y),z)=T^*(x,T^{*}(y,z))$ since $T^*$ is associative. Then
$$F_{M}(m_{i\oplus j,k})\le F_{M}(T^{*}(x,T^{*}(y,z)).$$
On the one hand, from Eqs.\eqref{eq4.8}, \eqref{eq4.9} and \eqref{eq4.13} we have
$$F_{M}(m_{i\oplus j,k})=c_{(i \oplus j)\oplus k}=c_{i \oplus j} \otimes c_k=(c_i\otimes c_j) \otimes c_k=(x\otimes y)\otimes z.$$
On the other hand, since $T^*(y,z) \in M$, we have $y \otimes z=T^*(y,z)$. Consequently,
$$F_{M}(T^*(x,T^*(y,z))=F_{M}(T^*(x,y\otimes z)=x\otimes(y\otimes z).$$
Therefore, $(x\otimes y)\otimes z \leq x\otimes (y\otimes z)$.

Subcase 3: The proof is completely analogous to Subcase 2 when $T^*(x,y) \in \cup_{l\in I}(a_{l},b_{l})$ and $T^*(y,z) \notin \cup_{l\in I}(a_{l},b_{l})$.
\end{proof}

In what follows, we discuss the case when $f(0^+)=0$. Write

(D8') $(\{c_i \mid i\in I\}, \otimes)$ is a semigroup where $\otimes: \{c_i \mid i\in I\} \times \{c_i \mid i\in I\} \to \{c_i \mid i\in I\}$ is the restriction of the operation $\otimes: M \times M \to M$.

Then the proof of the following two results are completely analogous to Propositions \ref{prop4.1} and \ref{propth4.1}, respectively.

\begin{proposition}\label{prop4.2}
 Let $(\mathcal{S}',C')=(\{[b_k', d_k']\mid k\in K\}, \{c_k' \mid k\in K\})$ be associated with $M \in \mathcal{A}_0$, $T^*:[0,1]^2 \to [0,1]$ be a continuous t-conorm and $f:[0,1] \to [0,1]$ be a strictly increasing function with $f(0^+)=0$ being an idempotent element of $T^*$. If there exist two elements $m_1,n_1 \in [0,1]$ with $m_1<n_1$ such that $f(x) \in (m_1,n_1)$ for each $x\in (0,1)$ and $\langle m_1,n_1,T^*_1 \rangle$ is the summand of the continuous t-conorm $T^*$, then (D2)-(D4) imply (D6), (D7) and (D8').
 \end{proposition}

\begin{proposition}\label{propth4.2}
Let $M \in \mathcal{A}_0$, $T^*:[0,1]^2 \to [0,1]$ be a continuous t-conorm and $f:[0,1] \to [0,1]$ be a strictly increasing function with $f(0^+)=0$ being an idempotent element of $T^*$. If there exist two elements $m_1,n_1 \in [0,1]$ with $m_1<n_1$ such that $f(x) \in (m_1,n_1)$ for each $x\in (0,1)$ and $\langle m_1,n_1,T^*_1 \rangle$ is the summand of the continuous t-conorm $T^*$, then $(M, \otimes)$ is a semigroup if and only if (D6), (D7) and (D8') hold.
\end{proposition}

The next theorem follows from Propositions \ref{prop3.1} and \ref{propth4.1} (resp. Proposition \ref{propth4.2}).

\begin{theorem}\label{Thm4.1}
Let $T^*:[0,1]^2 \to [0,1]$ be a continuous t-conorm, $f:[0,1] \to [0,1]$ be a strictly increasing function with $f(0^+) \notin \mbox{Acc}_{+}((\emph{Ran}(f))^{c})$ and $f(0^+)>0$ (resp. $f(0^+)=0$) being an idempotent element of $T^*$. If there exist two elements $m_1,n_1 \in [0,1]$ with $m_1<n_1$ such that $f(x) \in (m_1,n_1)$ for each $x\in (0,1)$ and $\langle m_1,n_1,T^*_1 \rangle$ is the summand of the continuous t-conorm $T^*$, then the generated function $T$ via Eq.\eqref{eq1.1} is a border continuous t-conorm if and only if (D6), (D7) and (D8) (resp, (D8')) hold.
\end{theorem}

Note that from the proof of Proposition \ref{prop4.1}, we know that (D7) is equivalent to Eqs.\eqref{eq4.5} and \eqref{eq4.6}. Hence, if (D6) holds in Proposition \ref{propth4.1} (resp. Proposition \ref{propth4.2}), then $(M,\otimes)$ is a semigroup if and only if  Eqs.\eqref{eq4.5} and \eqref{eq4.6} and (D8) (resp. (D8')) hold. Therefore, we have the following corollary.
\begin{corollary}
Let $M \in \mathcal{A}_0$ be determined by a sequence of triples $\{(a_i,b_i,c_i)\}_{i\in I}$, $T^*:[0,1]^2 \to [0,1]$ be a continuous t-conorm and $f:[0,1] \to [0,1]$ be a strictly increasing function with $f(0^+)>0$ (resp. $f(0^+)=0$) being an idempotent element of $T^*$. If there exist two elements $m_1,n_1 \in [0,1]$ with $m_1<n_1$ such that $f(x) \in (m_1,n_1)$ for each $x\in (0,1)$ and $\langle m_1,n_1,T^*_1 \rangle$ is the summand of the continuous t-conorm $T^*$, then $(M, \otimes)$ is a semigroup if and only if the following three conditions are satisfied:

$(i)$ $(\cup_{i,j \in I}(\min\{T^*(a_i,b_j),T^*(a_j,b_i)\},T^*(c_i,c_j)) \cap (\cup_{i \in I}(a_i,b_i))=\emptyset$;

$(ii)$ $(\cup_{i,j \in I}(T^*(a_i,a_j),T^*(b_i,b_j))) \cap \{a_i \mid i \in I\}=\emptyset$;

$(iii)$ $(\{c_i \mid i\in I\cup \{0\}\}, \otimes)$ is a semigroup where $\otimes: \{c_i \mid i\in I\cup \{0\}\} \times \{c_i \mid i\in I\cup \{0\}\} \to \{c_i \mid i\in I\cup \{0\}\}$ is the restriction of $\otimes: M \times M \to M$ (resp. $(\{c_i \mid i\in I\}, \otimes)$ is a semigroup where $\otimes: \{c_i \mid i\in I\} \times \{c_i \mid i\in I\} \to \{c_i \mid i\in I\}$ is the restriction of the operation $\otimes: M \times M \to M$).
\end{corollary}

The following example illustrates Theorem \ref{Thm4.1}.
\begin{example}\label{exp4.2}
\emph{$(1)$ Let the function $T^*:[0,1]^2\rightarrow [0,1]$ be defined by
\begin{equation*}
T^*(x,y)=\begin{cases}
(1-\frac{1}{2})S_p(\frac{x-\frac{1}{2}}{1-\frac{1}{2}},\frac{y-\frac{1}{2}}{1-\frac{1}{2}})+\frac{1}{2} & \text{if } (x,y)\in [\frac{1}{2},1]^2,\\
\max\{x,y\} & \text{otherwise}
\end{cases}
\end{equation*}
and the function $f:[0,1] \rightarrow [0,1]$ be defined by
\begin{equation*}
f(x)=\begin{cases}
\frac{1}{6}(x+3) & \text{if } x\in [0,1),\\
1 & \text{if } x=1,
\end{cases}
\end{equation*}
respectively.}

\emph{Then one can check that $M=[\frac{1}{2},\frac{2}{3}) \cup \{1\}$, $M$ is determined by the sequence of triples $\{(\frac{1}{2},\frac{2}{3},1)\}$, $f(0^+)=\frac{1}{2}$ is an idempotent element of $T^*$,  $f(0^+)$ satisfies the condition that $f(0^+) \notin \mbox{Acc}_{+}(M^c)$ and $f(x) \in (\frac{1}{2},1)$ for each $x\in (0,1)$. Meanwhile, $T^*=(\langle \frac{1}{2},1,S_p\rangle)$ and (D7) and (D8) hold. Therefore, from Theorem \ref{Thm4.1} we know that
\begin{equation*}
T(x,y)=\begin{cases}
x+y-\frac{1}{3}xy & \text{if } (x,y)\in [0,1)^2 \mbox{ and } x+y-\frac{1}{3}xy\in [0, 1),\\
1 & \text{otherwise}
\end{cases}
\end{equation*}
is a border continuous t-conorm.}

\emph{$(2)$ Let the function $S_L:[0,1]^2\rightarrow [0,1]$ be defined by
\begin{equation*}
S_L(x,y)=\min\{x+y,1\},
\end{equation*}}
\emph{the function $T^*:[0,1]^2\rightarrow [0,1]$ be defined by
\begin{equation*}
T^*(x,y)=\begin{cases}
\frac{1}{2}S_L(\frac{x-0}{\frac{1}{2}-0},\frac{y-0}{\frac{1}{2}-0}) & \text{if } (x,y)\in [0,\frac{1}{2}]^2,\\
\max\{x,y\} & \text{otherwise,}
\end{cases}
\end{equation*}
and the function $f:[0,1] \to [0,1]$ be defined by
\begin{equation*}
f(x)=\frac{1}{2}x,
\end{equation*}
respectively.}

\emph{Then one can check that $M=[0,\frac{1}{2}]$, $M$ is determined by the sequence of triples $\{(0,\frac{1}{2},\frac{1}{2})\}$, $f(0^+)=0$ is an idempotent element of $T^*$ and $f(0^+)$ satisfies that $f(0^+) \notin \mbox{Acc}_{+}(M^c)$, $f(x) \in (0,\frac{1}{2})$ for each $x\in (0,1)$. Meanwhile, $T^*=(\langle 0,\frac{1}{2},S_L\rangle)$ and (D7) and (D8') are satisfied. Therefore, from Theorem \ref{Thm4.1} we know that
\begin{equation*}
T(x,y)=\min\{x+y,1\}=S_L(x,y)
\end{equation*}
is a border continuous t-conorm.}
\end{example}

\section{Conclusions}

The main contribution of this article is that we supplied a necessary and sufficient condition for the generated function $T$ via Eq.\eqref{eq1.1} being a border continuous t-conorm when $f:[0,1] \to [0,1]$ is a strictly increasing function (Theorem \ref{Thm4.1}). Dually, one can give a necessary and sufficient condition for the generated function $T$ via Eq.\eqref{eq1.1} being a border continuous t-norm when $f:[0,1] \to [0,1]$ is a strictly decreasing function. Therefore, we completely characterize the border continuous triangular conorm $T$ via Eq.\eqref{eq1.1} when $f$ is a strictly monotone function. Comparing with Vicen\'{\i}k \cite{PV2008}, the different in essence is that his work holds over the semigroup $([0,\infty],+)$ with the usual addition operator ``+", and our results are true when $T^*$ can be written as an ordinal sum of continuous Archimedean t-conorms.
%\section*{Acknowledgments}
%The authors are grateful to the anonymous referees for their valuable comments and suggestions, which helps the authors to improve the earlier version of this article.

\section*{Compliance with Ethical Standards}
\quad\\
{\bf{Ethical approval}} This article does not contain any studies with human
participants or animals performed by any of the authors.\\
{\bf{Funding}} This work was supported by the National Natural Science Foundation of China under Grant No.12471440.\\
{\bf{Conflict of interest}} All the authors in the paper have no conflict of
interest.\\
{\bf{Informed consent}} Informed consent was obtained from all individual
participants included in the study.

\end{document}